\documentclass[10pt,twocolumn,letterpaper]{article}

\usepackage[pagenumbers]{cvpr} % To force page numbers, e.g. for an arXiv version

\usepackage{graphicx}
\usepackage{amsmath}
\usepackage{amssymb}
\usepackage{booktabs}
\usepackage{array}
\usepackage{amsmath}

\usepackage{cuted}
\usepackage{capt-of}
\usepackage{dsfont}
\usepackage{listings}
\usepackage{pythonhighlight}

\newcommand{\PreserveBackslash}[1]{\let\temp=\\#1\let\\=\temp}

\newcolumntype{C}[1]{>{\PreserveBackslash\centering}p{#1}}
\newcolumntype{R}[1]{>{\PreserveBackslash\raggedleft}p{#1}}
\newcolumntype{L}[1]{>{\PreserveBackslash\raggedright}p{#1}}

\newcommand\scalemath[2]{\scalebox{#1}{\mbox{\ensuremath{\displaystyle #2}}}}

\usepackage[pagebackref,breaklinks,colorlinks]{hyperref}

\usepackage[capitalize]{cleveref}
\crefname{section}{Sec.}{Secs.}
\Crefname{section}{Section}{Sections}
\Crefname{table}{Table}{Tables}
\crefname{table}{Tab.}{Tabs.}

\begin{document}

%%%%%%%%% TITLE - PLEASE UPDATE
\title{Analysis of Hospital Bed Requirements Using Discrete Event Simulation and Mathematical Modeling}

\author{Dincer Atasoy\\
Igdir University\\
Igdir, Turkey\\
{\tt\small dincer.atasoy@igdir.edu.tr}
}

\maketitle

%%%%%%%%% ABSTRACT

\begin{abstract}
Using SimPy and Discrete Event Simulation we have observed the different model responses of a system consisting of a hospital and people getting sick/healing under different initial conditions.

In our model, each independent person can get sick at an exponential rate.  The hospital's capacity is limited and accepts only a limited amount of sick people. When the hospital gets full, people are sent home to heal. This is simulated under different conditions and the simulation results have demonstrated that the model reaches the steady-state in every scenario. These cause the outputs of the simulation to be similar to each other even when the simulation is run with different initial conditions. The results are compared with the Machine Repair Problem and Erlang-Loss System and the model is validated.
\end{abstract}

\section{Introduction}
Analyzing large real systems is complex and it is impossible to observe the system under different conditions. Discrete Event Simulation\cite{fishman2013discrete} helps people understand the behavior and mechanism of a real system by simulating it under certain conditions. The simulation displays each step of the model and creates conclusions about the system. In a simulation entities and processes interact with each other, causing different results. This way, the real-time behavior of the system can be modeled as events occur. 

This study simulates a real-time system consisting of a hospital and people getting sick/healing. The simulation demonstrates different outputs of the model as the hospital status/simulation time changes. In this study,  the model outputs and responses are measured and discussed. The problem will be stated and the model will be explained in the next chapter. Following that, the numerical results will be shared and discussed in the scope of Queueing Systems\cite{white2012analysis} more precisely Machine Repair Problem\cite{posafalvi1989numerical} and Erlang-Loss System\cite{erlang1909theory}. Finally, the conclusion of the study will be stated.

\section{Related Work}
\label{sec:related}
Discrete event simulation is widely used for a lot of reasons and allows users to better understand its behavior and mechanism by replicating the behavior and mechanism of a genuine system under specific conditions. Unlike continuous event simulation, it suits making computations efficiently, and there are works that compare discrete event simulation with continuous event simulation in certain aspects \cite{doi:10.1080/17477778.2018.1465153}.

Until now, there has been research \cite{doi:10.1080/17477778.2019.1664264} about using discrete event simulation to understand the underlying structures of how diseases are spreading. However, the amount of research has dramatically increased in recent years due to the COVID-19 virus.

Recent work \cite{bartzbeielstein2020hospital} introduces babsim.hospital, a resource-planning tool for hospitals dealing with the COVID-19 epidemic. Another work \cite{healthcare10020189} intends to assess the effects of modifying consultation start time and patient arrival on wait times and congestion in a dual practice outpatient clinic. At the same time, a work \cite{10.1371/journal.pone.0253869} seeks to estimate the appropriate number of machines and operators needed, as well as their placement at various workstations, based on the available resources and the rate of samples to be examined every day. Finally, the goal of recent work \cite{9384093} is to propose a decision-aid tool for hospital management that will allow them to decide on the bed needs for a specific hospital or network of hospitals over a short-to-medium term horizon.

\section{Methodology}

We aim to model a system where there are $N=1582$ people that each can get sick with an exponential rate of $\lambda = 1/300\ [patients/day]$. People who are sick can go to the hospital to heal or they can stay home. The probability of a sick person going to the hospital for healing is $0.2$ , whereas the probability of a sick person staying home is $0.8$.

Our simulation is a process-based simulation where processes interact with each other as objects with their own states and variables. There is an environment under which each process occurs i.e. people get sick/heal.  And processes communicate through $yield$ and wait for each other. Also, there are several processes active throughout the simulation. There is a simulation i.e. environment time running at all times, waiting for the events. 

 If our model was to be simulated as an event-based simulation, then there would be a list of events and the simulation would change its states as events occur according to their arrival time. In this case, we would have to implement an event scheduler, so that when people arrive, the system state changes. Since there is an event queue in an event-based case, only one routine would be active at a time. But in our case, there are several processes active in a time period, they work in threads. Also, in our model, we check the availability of the hospital by using the communication of processes and objects. We can easily do this with SimPy\cite{matloff2008introduction}. In an event-based approach, there would be no need for SimPy since there are no processes interacting with each other. So, we would use normal Python loops and variables to control the number of people in the hospital and decide whether the system should accept patients to the hospital or not by occurring each event. In event-based, keeping the FEL would be easier, since each event is listed in the beginning. 

In both approaches, the results won't change, only the implementation would change. So, one can decide on which simulation type to use by deciding the desired implementation.

There are $K=N/24=66$ beds in the hospital which means that only 66 people can be healed at the hospital. After the hospital is full, the patients arriving at the hospital start to get rejected. Those patients are sent home to heal. Since the healing process at the hospital takes less time, the healing duration of a sick person at home takes $r \sim U[1, 2]$  times longer. The unit of time is days.

There are 3 different exponential healing rates for 3 different types of patients: 

\begin{enumerate}
\item $\mu _1:$  Healing rate of the patient being healed at the hospital. 
\item $\mu_2:$  Healing rate of the patient healing at home. 
\item $\mu_3:$  Healing rate of the patient that is initially meant to be healed at the hospital but sent home due to full capacity. So s/he is obligated to heal at home. 
\end{enumerate}

The healing rates are:

\begin{equation}
  \begin{array}{l}
\mu_1=1/6\\
\mu_2=1/10\\
\mu_3^{-1}=\mu_1^{-1}\cdot r=6\cdot r\text{ where }r \sim U[1,2]
  \end{array}
\end{equation}

We have a system as described above and we are to generate a Discrete Event Simulation for it. We have used SimPy\cite{matloff2008introduction} to simulate the events and observe the results. Events can be getting sick and a healing process is done for each person in the system. Everyone in the system can get sick and once a person is healed s/he can get sick again. For the simulation to be consistent with reality, we choose people who get sick randomly. So, anyone can get sick in the system, no one is favored.

The simulation runs in a SimPy environment called $env$. In that environment, the resource is the hospital with a capacity of 66, since each sick person arriving at the hospital expects to be healed there. If the resource is occupied, i.e. if there are 66 people in the hospital, then next arrivals are sent home, and until someone leaves the hospital sick people coming to the hospital get rejected. The simulation run time is 1000, 10000, or 100000 units of time. 

So there are two subsystems in our model: hospitals and people's homes.

To start to run the simulation, one should set the $until$ and $hospital\_status$  to the desired simulation run time and hospital status. This can be done by changing the commented out variables at the beginning of the code:

\begin{python}
# until=1000 
until=10000 
# until=100000 

# hospital_status = "empty" # state 
# hospital_status = "half-full" # state 
hospital_status = "full" # stat

# RANDOM_SEED = 978
RANDOM_SEED = 979
\end{python}

Then, run the simulation by running the code in .ipynb file. 

Firstly, $1582$  people are created for the simulation. Each person lives in the same environment, has unique IDs and information on where to heal, heal duration, and departure time. 

Secondly, the hospital is filled according to the $hospital\_status$.The system adds the necessary number of people to the hospital before the simulation starts to run. 

Finally, the simulation starts in the environment by calling the generator of sick people ($person\_generator()$) and runs until the simulation time is done. During the runtime of the simulation, several people get sick, get healed and get sick again, and so on.  The simulation starts after one interarrival time. 

The $person\_generator()$, first checks whether the system is full of sick people or not. If everyone in the system is sick, then it waits for one to heal by yielding the $generator\_reactivate$ to reactivate the system. It then chooses a random healthy person to get sick (again) using the list $healthy\_people$. This person is removed from the list of healthy people and his/her interarrival time is calculated using ${\lambda\cdot |healthy\ people|}$ as an exponential distribution\cite{marshall1967multivariate}. After that, this person is added to the $will\_arrive\_people$ list including his/her arrival time and id. 

The simulation now waits for the selected sick person to arrive, i.e. it waits for his/her arrival time. It removes them from the healthy people list and Future Event List. At this step, the simulation decides whether that person will go to the hospital or not. After that, this person is processed to be sick and sent to $get\_sick()$ method.

In $get\_sick()$, the healing place of the patient is checked. This is the part of the model where a person gets sick and decides to go to the hospital or not. 

If he/she decided to go to the hospital, then the capacity of the hospital is checked. This is done by checking the $hospital\_count$ where the resource capacity is checked. 

If the hospital is not full, then the patient gets accepted to the hospital. A hospital bed is requested ($hospital.request()$) and yielded. The heal duration is calculated using $\mu_1,$ $random.expovariate(mu_1)$.

Else if the hospital is full, the person is sent home to heal. The heal duration is calculated using $\mu_3$ .

If he/she decides to heal at home, then they start healing at home. The heal duration is calculated using $\mu_2,$ $random.expovariate(mu_2)$.

Each sick person is then prepared for healing. In $prep\_to\_heal()$ the departure time of the person is calculated and the sick person is added to the sick people list. The data of that person is added to the table along with the FEL. 

Finally, the person starts to heal. In $heal()$ , the simulation waits for the heal duration of the person. At this point, if everyone in the simulation is sick, then the simulation reactivates the hospital since one patient will be leaving. That person is removed from the sick list and added to the healthy list. The data of the recovered person is added to the table.

In our simulation, the output is an Excel table containing the important DES values as events occur. It contains Simulation Time, Person Id, Event Type, Healing Place, Heal Time, Number of Sick People, Number in the Hospital, Future Event List. These Excel tables can be found in Section \ref{sec:num}. 

\begin{figure*}
\begin{center}
\centerline{\includegraphics[width=2\columnwidth]{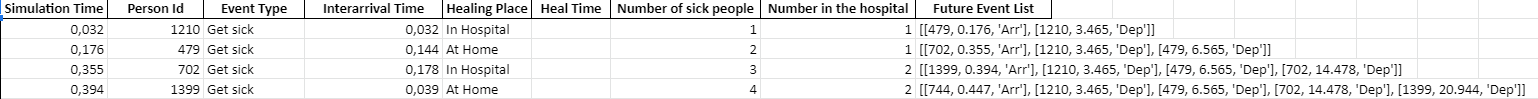}}
\caption{An example of the DES Table in Excel form (Seed = 978, Hospital=empty)}
\label{fig:1}
\end{center}

\end{figure*}

\begin{figure}
\begin{center}
\centerline{\includegraphics[width=1\columnwidth]{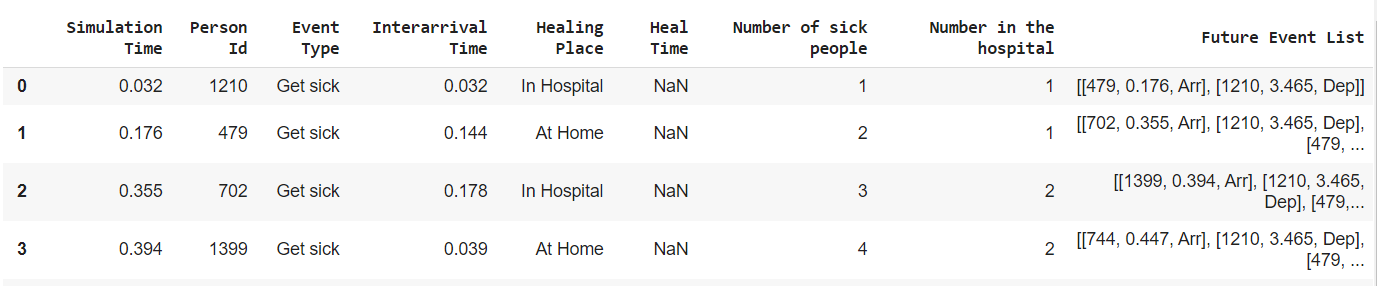}}
\caption{An example of the DES Table in the code (Seed=978, Hospital=empty)}
\label{fig:2}
\end{center}

\end{figure}

\textbf{Explanation of the DES Table:}

\begin{enumerate}
\item Person with id 1210 gets sick and he/she is to be healed in the hospital. The number of sick people and the number in the hospital increase by 1. The next arrival and the current person's departure time are seen in FEL. The departure time = simulation time + healing duration. Also, the interarrival time for the current sick person is written in the table in Figure \ref{fig:1}.
\item Person with id 479 arrives and will be healed at home. The number of sick people increases by 1 and the number in the hospital remains unchanged. The next arrival and the current person's departure time are added to the FEL in Figure \ref{fig:2}.
\item In these examples, Heal Time is NaN since that value is written in the table only for healing events, i.e. when a sick person departs.
\end{enumerate}

The simulation also calculates the following responses:

\begin{itemize}
\item long-run probability of the hospital is empty,
\item proportion of sick people healing in the hospital,
\item average number of sick people in the population,
\item average proportion of sick people in the population,
\item standard deviation of sick people in the population,
\item average number of occupied beds in the hospital,
\item standard deviation of occupied beds in the hospital
\item total average sickness time.
\end{itemize}

The simulation also plots the change in the number of sick people in the system and the number of people in the hospital over Simulation Time (See \ref{analysis:1} and \ref{analysis:2}).   

The analysis of the outputs, graphs, and responses can be found in Section \ref{sec:num}.

\begin{figure}
\begin{center}
\centerline{\includegraphics[width=1\columnwidth]{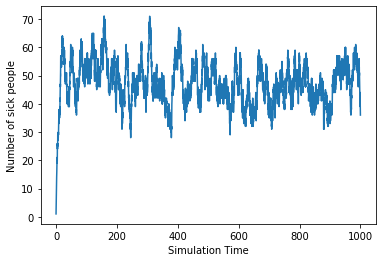}}
\caption{Hospital Status = Empty, Number of sick people}
\label{fig:3}
\end{center}

\end{figure}

\begin{figure}
\begin{center}
\centerline{\includegraphics[width=1\columnwidth]{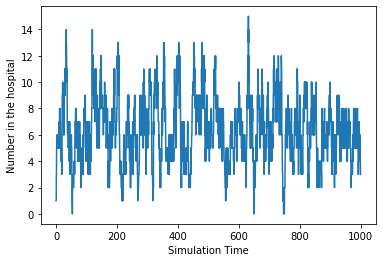}}
\caption{Hospital Status = Empty, Number in the hospital}
\label{fig:4}
\end{center}

\end{figure}

\begin{figure}
\begin{center}
\centerline{\includegraphics[width=1\columnwidth]{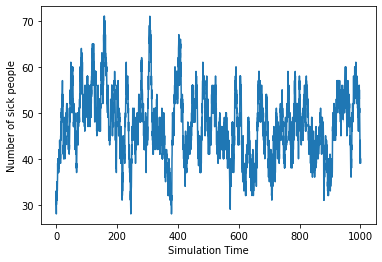}}
\caption{Hospital Status = Half Full, Number of sick people}
\label{fig:5}
\end{center}

\end{figure}

\begin{figure}
\begin{center}
\centerline{\includegraphics[width=1\columnwidth]{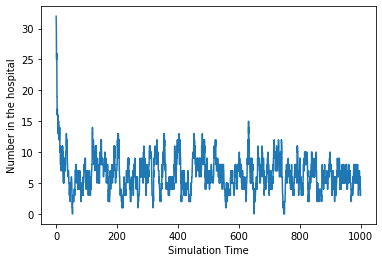}}
\caption{Hospital Status = Half Full, Number in the hospital}
\label{fig:6}
\end{center}

\end{figure}

\begin{figure}
\begin{center}
\centerline{\includegraphics[width=1\columnwidth]{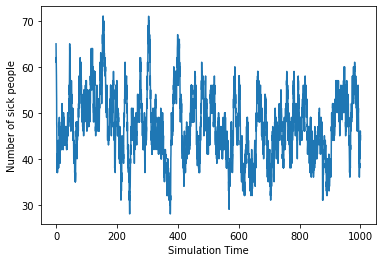}}
\caption{Hospital Status = Full, Number of sick people}
\label{fig:7}
\end{center}

\end{figure}

\begin{figure}
\begin{center}
\centerline{\includegraphics[width=1\columnwidth]{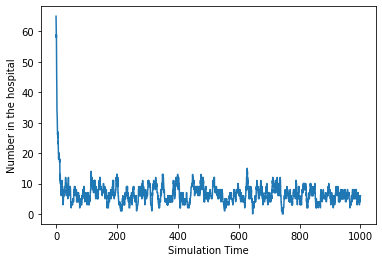}}
\caption{Hospital Status = Full, Number in the hospital}
\label{fig:8}
\end{center}

\end{figure}

\section{Numerical Analysis}
\label{sec:num}

\subsection{Analysis of the Plots}
\label{analysis:1}

As one can see from the graphs, with both seeds and each simulation run time-hospital status pair, the steady-state mean values are approximately the same. However, the standard deviation of the values gets larger as simulation run time gets larger which is the expected behavior of the system.  

Moreover, when the system starts with half-full/full hospital, the number of people in the hospital decreases at a fast rate and reaches a steady state. This fast fall in the number of people in the hospital is due to the fast healing rates of people in the hospital. And when the steady-state is reached, due to the fast healing of people in the hospital, the hospital doesn't reach its full capacity again. When people come to the hospital, they heal faster than the people healing at home ($\text{See }\mu_1, \mu_2$ ).  These and other responses will be validated and analyzed.

\begin{figure*}
\begin{center}
\centerline{\includegraphics[width=2\columnwidth]{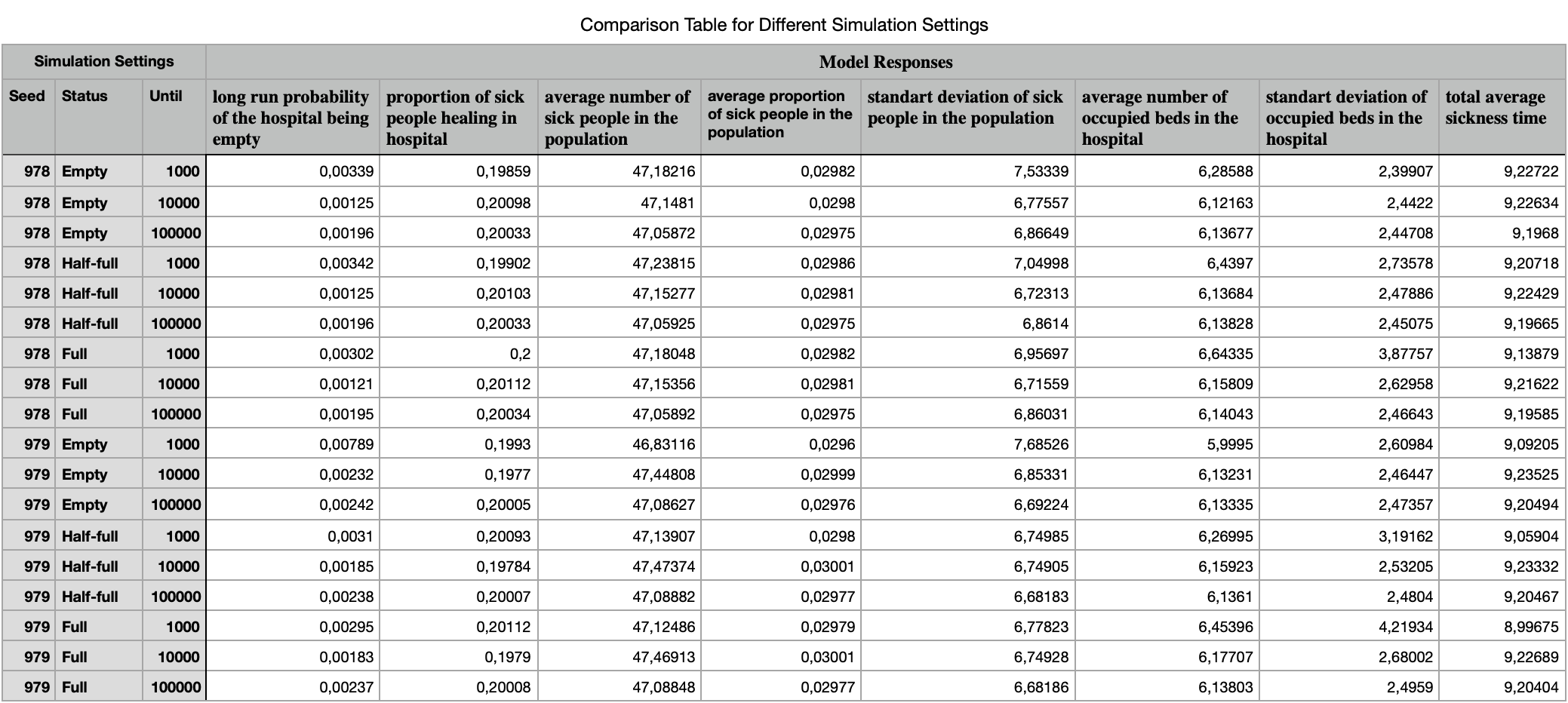}}
\caption{Comparison Table for Different Simulation Settings}
\label{fig:9}
\end{center}

\end{figure*}

\begin{figure}
\begin{center}
\centerline{\includegraphics[width=1\columnwidth]{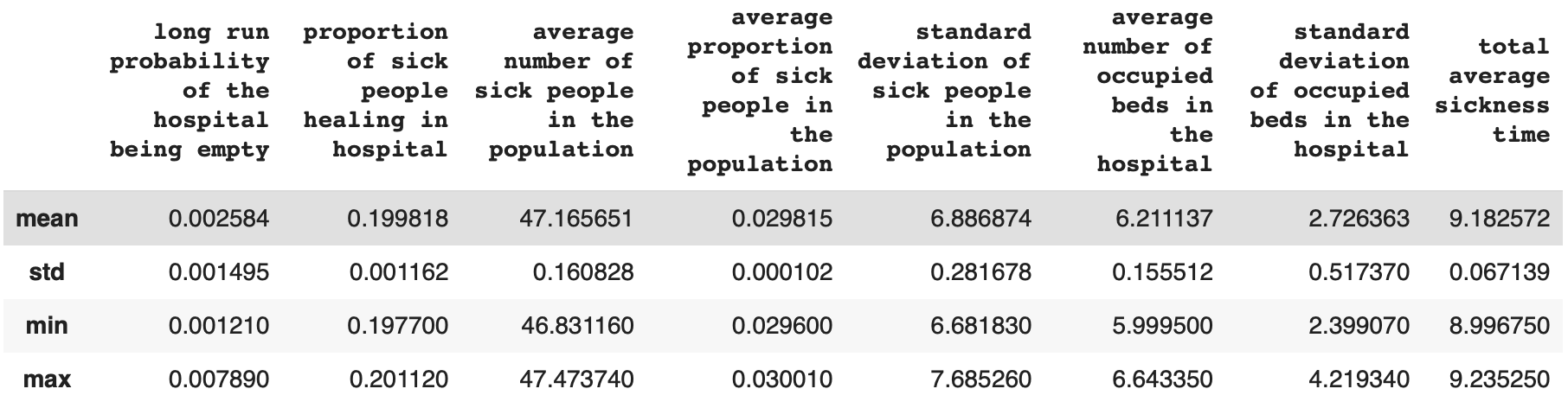}}
\caption{Summary Table for Different Simulation Settings}
\label{fig:10}
\end{center}

\end{figure}

\subsection{Validation of Model with help of Queueing Systems}
\label{analysis:2}

We can validate our model with Queueing Systems \cite{allen2014probability}. Our model could be an approximation to Machine Repair Problem \cite{sztrik1988g} because the arrival rate is changing with a number of healthy people throughout the simulation. Also, population size is limited. Then, model will be $M/M/N/\infty/N$ \cite{sztrik1985finite}

Also, we can model our hospital system as another model. Because it has K beds and the person who couldn't get help from the hospital will not wait for the hospital, he/she will be healed at their home. Thus, the system will be $M/M/k/k\text{ Queue(Erlang-Loss System)}$ \cite{erlang1909theory}.

Let's first calculate the common values of both systems:

Let $\lambda = 1/300$  , $N=1582$

There are three different $\mu$'s in our case, so we have to find the expected value of $\mu$ in order to validate Machine Repair Problem responses. However, the probability of occurring $\mu_1,\mu_3$ depends on the probability of a person being rejected by the hospital, and this probability is dependent on$\lambda_e$ (because the arrival rate of the hospital is $\lambda_e\cdot0.2$).  Therefore, we will look at the edge cases of the probability to calculate an interval for it.
If people who need to go to the hospital get accepted by the hospital i.e. capacity is not full:

\begin{equation}
  \scalemath{0.8} {
  \begin{array}{l}
\mu=\frac{1}{E[S]}=\frac{1}{E[S_1]\cdot0.2+E[S_2]\cdot0.8 + E[S_3]\cdot0} = \frac{1}{6\cdot0.2+10\cdot0.8} = 0.109
  \end{array}
  }
\end{equation}

If people who need to go to the hospital get rejected by the hospital i.e. capacity is full:
        
Let $R \in U[1,2]$, then since $\mu_1$ and R are independent:
\begin{equation}
  \begin{array}{l}
E[S_3]=E[S_1]\cdot E[R]=61.5=9
  \end{array}
\end{equation}

\begin{equation}
\scalemath{0.8} {
  \begin{array}{l}
\mu=\frac{1}{E[S]}=\frac{1}{E[S_1]\cdot0+E[S_2]\cdot0.8 + E[S_3]\cdot0.2}= \frac{1}{9\cdot0.2+10\cdot0.8} = 0.102
  \end{array}
}
\end{equation}

Then, $0.102 \leq \mu \leq 0.113$ so we can continue to calculate values for two edge cases of service rate which are 0.102 and 0.113
    
Case 1 $: \mu = 0.102$

\begin{equation}
\scalemath{0.7} {
  \begin{array}{l}
        {P_0=\Biggr [
        \sum\limits_{n=0}^{c-1}{N \choose n}\cdot {\Bigr (\frac{\lambda}{\mu}\Bigr)}^n+ \sum\limits_{n=c}^{N}\Bigr (\frac 1 {c!}\Bigr)\cdot{N \choose n}\cdot(\frac {n!} {c^{(n-c)}})\cdot{\Bigr (\frac{\lambda}{\mu}\Bigr)}^n\Biggr]^{-1}} \\
        
        \ \\
        
        P_0=\Biggr [
        \Bigr (\frac 1 {1582!}\Bigr)\cdot{1582 \choose 1582}\cdot(\frac {1582!} {1582^{(1582-1582)}})\cdot{\Bigr (\frac{1/300}{0.102}\Bigr)}^{1582}
         +
        \sum\limits_{n=0}^{1581}{1582 \choose n}\cdot {\Bigr (\frac{1/300}{0.102}\Bigr)}^n \Biggr]^{-1} \\
        
        \ \\
        
        \implies P_0= 8.06\cdot 10^{-23} \\
        
        \ \\
        
        P_n=\frac{N!/(N-n)!}{n!}(\lambda / \mu)^n\cdot P_0 \\ \ \\
        \implies P_n=\frac{1582!/(1582-n)!}{n!}((1/300) / 0.102)^n\cdot 8.06\cdot 10^{-23} \\
        
        \ \\
        
        L_{0.102}=\sum \limits_{n=0}^Nn\cdot P_n=50.06 \\
        
        \ \\
        
        \implies \lambda_e=\lambda(N - L)=\frac{1}{300}\cdot (1582-50.06)=5.106 \\
  \end{array}
  }
\end{equation}

Case 2$: \mu=0.109$

\begin{equation}
  \scalemath{0.7} {
  \begin{array}{l}

        {P_0=
        \Biggr [\sum\limits_{n=0}^{c-1}{N \choose n}\cdot {(\frac{\lambda}{\mu})}^n+ \sum\limits_{n=c}^{N}(\frac 1 {c!})\cdot{N \choose n}\cdot (\frac {n!} {c^{(n-c)}})\cdot{(\frac{\lambda}{\mu})}^n\Biggr ]^{-1}} \\
        
        \ \\
        
        P_0=\Biggr [\sum\limits_{n=0}^{1581}{1582 \choose n}\cdot {(\frac{(1/300)}{0.109})}^n+
        \sum\limits_{n=1582}^{1582}(\frac 1 {1582!})\cdot {1582 \choose n}\cdot \Big(\frac {n!} {1582^{(n-1582)}}\Big)\cdot {\Big(\frac{\lambda}{\mu}\Big)}^n\Biggr ]^{-1}\\
        
        \ \\
        
        {\implies P_0=\Biggr [\Big(\frac{1/300}{0.109}\Big)^{1582}+\sum\limits_{n=0}^{1581}{1582 \choose n}\cdot \Big (\frac{1/300}{0.109}\Big)^n\Biggr ]^{-1}=2.01\cdot10^{-21}} \\
        \ \\
        
        P_n=\frac{1582!/(1582-n)!}{n!}\Big (\frac{1/300} {0.109}\Big)^n\cdot2.01\cdot10^{-21} \\
        
        \ \\
        
        L_{0.109}=\sum \limits_{n=0}^{1582}n\cdot P_n=46.85 \\
        
        \ \\
        
        \implies \lambda_e=\lambda(N - L)=\frac{1}{300}\cdot(1582-46.85)=5.118 \\
  \end{array}
  }
\end{equation}

Let's now validate the hospital model:

Let $\lambda_h$ be the arrival rate to the hospital. Then, 

for $\lambda_{e,0.102}=5.106$

\begin{equation}
  \begin{array}{l}
    \lambda _ h=\lambda_{e,0.102} \cdot (0.2)=\lambda_{e,0.102} /5=\frac {5.106}{5}\approx1.02
  \end{array}
\end{equation}

for $\lambda_{e, 0.109}=5.118$

\begin{equation}
  \begin{array}{l}
  \lambda_h=\lambda_{e,0.109}\cdot (0.2)=\lambda_{e,0.109}/5=\frac{5.118}{5}\approx1.02
  \end{array}
\end{equation}    

\begin{enumerate}
\item Long run probability of hospital being empty $(P_{0,h})$:

\begin{equation}
  \scalemath{0.8} {
  \begin{array}{l}
    P_{0,h}=\frac{1}{\sum\limits_{k=0}^{n} ((\lambda_h/\mu_1)^k/k!)}=\frac{1}{\sum\limits_{k=0}^{66} ((\frac{1.02}{1/6})^k/k!)}=0.0022
  \end{array}
  }
\end{equation}

    This matches the simulation result for the long-run probability of the hospital being empty which was 0.0025. Then, the response is validated.

\item Average number of beds occupied in the hospital $(L_h)$:

\begin{equation}
  \begin{array}{l}
  W_h=\frac{1}{\mu_1}=6 \\ \ \\
    L_h=\lambda_h\cdot W_h=1.02\cdot6=6.12
  \end{array}
\end{equation}
    
    This matches the simulation result for the average number of occupied beds in the hospital which was 6.21. 

\item Possibility of people who go to the hospital being rejected $(P_{66})$:

\begin{equation}
  \scalemath{0.8} {
  \begin{array}{l}
  P_k=P_{0,h}\cdot(\lambda_h/\mu_1)^k\cdot \frac 1 {k!} \\ \ \\
    P_{66}=0.0022 \cdot\Big(\frac{1.02}{1/6}\Big)^{66}\cdot \frac 1 {66!}=3.405\cdot10^{-44}\approx0
  \end{array}
  }
\end{equation}

    Then, the possibility of people who go to the hospital being rejected is 0. Therefore, $\mu_3$ is not important in our system, it can be omitted.
\end{enumerate}

In conclusion, $\mu$ is taken as

\begin{equation}
  \scalemath{0.8} {
  \begin{array}{l}
 \mu=\frac{1}{E[S]}=\frac{1}{E[S_1]\cdot0.2+E[S_2]\cdot0.8 + E[S_3]\cdot0}= \frac{1}{6\cdot0.2+10\cdot0.8} = 0.109 
  \end{array}
  }
\end{equation}

and the effective arrival rate is taken as:

\begin{equation}
  \begin{array}{l}
\lambda_e=\lambda(N - L)=\frac{1}{300}\cdot(1582-46.85)=5.118
  \end{array}
\end{equation}

Let's now validate the rest of the model using the Machine Repair Problem:

The theoretical/calculated values are compared with Figure \ref{fig:9}.

\begin{enumerate}

\item 
The average number of sick people in the population $(L):$
    
    We have calculated above that:

\begin{equation}
  \begin{array}{l}
    L=L_{0.109}=\sum \limits_{n=0}^{1582}n\cdot P_n=46.85

  \end{array}
\end{equation}

    This matches the result of the simulation where the mean value of L was 47.16.
    
    Hence, our model's average number of sick people is consistent with the theoretical value. 

\item
Total average sickness time $(W):$

\begin{equation}
  \begin{array}{l}
    W=\frac{L}{\lambda_e}=\frac{46.85}{5.122}=9.15

  \end{array}
\end{equation}

    This matches the result of the simulation for the total average sickness time which was 9.18.
    
    Hence, our model's total average sickness time is consistent with the theoretical value.

\item
The average proportion of sick people in the population $(\frac LN):$
    
    Average proportion of sick people in the population = $\frac{L}{N}= \frac{46.85}{1582} = 0.0296$
    
    This matches the result of the simulation for the average proportion of sick people in the population which was 0.0298.
    
    Hence, our model's average proportion of sick people in the population is consistent with the theoretical value.

\end{enumerate}

\subsection{Analysis of Model Responses}
\label{analysis:3}

We combined all data gathered by different settings in Figure \ref{fig:9} to analyze it, and we will analyze it one by one according to different settings. 

\begin{enumerate}

\item
long-run probability of the hospital being empty:
    - When we look at the values $P_{0,h}$ of different settings, all values are small probabilities, and they have a mean $0.002584$ (Figure \ref{fig:10}). It shows us that in the long run hospital is converging to be empty.
        
        The main reason for that is the arrival rate is much bigger than the service rate. Thus, the hospital is mostly not empty. $\lambda_e=5.118, \mu=0.109, L=46.85$ (calculated in \ref{analysis:2}). 
        
        Also, we can see the same from Figure \ref{fig:10} that the mean value of the average number of sick people in the population is $47.16$. So, it is a low possibility that the system will get close to being empty and so hospital as well.
        
    The standard deviation of the probability is large when compared to the mean value (almost the same). That means the variance of the probability is big.
    The main reason for that is since the probability value is low, the fluctuation of probability is big in small simulation time. (It is clear in above line plots (\ref{fig:3} and \ref{fig:4}.) which have simulation time equal to 1000 time units.)

\item
the proportion of sick people healing in the hospital:
    
    That response is straightforward to analyze. It has a small standard deviation and is always close to the mean value $0.2$. It is expected since the probability of a sick person getting treated in the hospital is given as $0.2$. Also, they have the same values because, in all settings, the number of beds in the hospital is not getting full. Therefore, a person who needs to heal in the hospital always has a chance to heal in the hospital. Hence, the final value of response is exactly the same as the defining value of needing treatment in the hospital.

\item
the average number of sick people in the population:
    
    We have validated the response above in Section \ref{analysis:2} The response is not fluctuating so much because the simulation never exceeds the total number of beds in the hospital, and for every setting, it is between $47$ and $48$, also has a very small standard variation. As a result, because of the service rate and arrival rate and also the number of beds in the hospital, the simulation is reaching its steady-state quickly in every setting. Besides, the simulation run-time is big enough for the simulation to reach its steady-state.

\item
the average proportion of sick people in the population:
    
    We have validated the response above in Section \ref{analysis:2} and it has exactly the same explanation with the average number of sick people in the population.

\item
the average number of occupied beds in the hospital
    
    We have validated the response above in Section \ref{analysis:2}. The response is not changing so much with different settings Because of the reason we mentioned above.

\item
total average sickness time:
    
    We have validated the response above in Section \ref{analysis:2}. The total average sickness time is equal to service time in our model, and it is generated with exponential distributions $\mu_1,\mu_2,\mu_3$ . However, as we analyzed in Section \ref{analysis:2}, $\mu_3$ is never the case in our model, and we can find the weighted average of service time. Moreover, when we take the multiplicative inverse of the service time, we get average sickness time and as we validated in Section \ref{analysis:2}, it is very close to the theoretical value.

\end{enumerate}

\subsection{Conclusion of Numerical Analysis}
\label{analysis:4}

The simulation is not changing in different settings. Because of four reasons:

\begin{enumerate}
\item
Simulation time is long enough to get our simulation to steady-state, even 1000 unit time is enough.
\item
The total number of beds in the hospital (K) is big enough for every person who needs to heal in the hospital to be able to heal in the hospital.
\item
Both arrival rate and service rate are small compared to the total simulation time. Therefore, we couldn't observe different responses in different simulation time settings.
\item
The ratio of arrival rate and service rate is not big enough to reach the total population size. Therefore, total population size is not a constraint for the simulation.
\end{enumerate}

Even though the simulation settings were not so great to observe different simulation responses in different settings, it is validated nicely with help of the Queueing Systems \cite{adan2002queueing}.

\section{Conclusion}
\label{sec:conclusion}

The study's goal was to create a model to observe the real-time behavior of a system in which each person can get sick and go to the hospital/stay home for healing. This model is simulated under different starting conditions and it is observed that overall results don't fluctuate from each other. The system reaches its steady-state in each simulation time and it is impossible for the system to fill the hospital during the run time since the ratio of arrival rates and service rates are relatively small compared to the number of beds in the hospital. Under a smaller hospital capacity or smaller simulation run time or higher arrival rates and lower healing rates, the results of the study would be different. Under these conditions, the model responses in each initial hospital status would also look more different than in our study. In the future, one can model a similar system with different initial conditions to get a variety of model responses.
 
%%%%%%%%% REFERENCES
{\small
\bibliographystyle{ieee_fullname}
\bibliography{main}
}

\appendix

\end{document}